\documentclass[%
 aip,
cp,  
 amsmath,amssymb,
 reprint,%
]{revtex4-2}

\usepackage{graphicx}
\usepackage{dcolumn}
\usepackage{bm}

\usepackage[utf8]{inputenc}
\usepackage[T1]{fontenc}
\usepackage{mathptmx} 

\begin{document}

\title{Numerical Modeling of Kondratyev's Long Waves Taking into Account Heredity}

\author{Makarov Danil} 
\email[Corresponding author: ]{danil.makarov.pk@yandex.ru}
\affiliation{Vitus Bering Kamchatka State University, Petropavlovsk-Kamchatskiy, Russia
}
\affiliation{Kamchatka Technical University, Petropavlovsk-Kamchatskiy, Russia}

\author{Parovik Roman}%
\email{romanparovik@gmail.com}
\affiliation{Institute of Cosmophysical Research and Radio Wave Propagation, Far Eastern Branch of the Russian Academy of Sciences, Kamchatka Territory, Paratunka, Russia
}

\date{\today} 

\begin{abstract}
The paper proposes a new mathematical model of economic cycles and crises, which generalizes the well-known model of Dubovsky S.V. The novelty of the proposed model lies in taking into account the effect of heredity (memory), as well as the introduction of harmonic functions responsible for the arrival of investments in fixed assets and new management technologies in innovation. The mathematical description is given using the Gerasimov-Caputo fractional derivatives, which are studied within the framework of the theory of fractional calculus. The mathematical model was investigated using the numerical method of Adams-Bashforth-Moulton (ABM), phase trajectories were constructed. It is shown that the proposed mathematical model can have both regular and chaotic regimes.
\end{abstract}

\maketitle

\section{Introduction}

In \cite{Tarasov_2019}, a literature review was carried out, about 260 sources, devoted to the application of fractional calculus in economics. It also reviewed the work of various authors who used fractional calculus to describe differential models of the economy, taking into account heredity or memory. Memory effects in an economic system manifest themselves in such a way that its current state depends on previous states or prehistory. Such effects can be taken into account when studying economic crises and cycles that arise under certain conditions, depending on the prehistory. Long waves of Kondratyev are not an exception \cite{Alexander_2002}. They describe well the innovation systems that characterize the explosive technological growth (breakthrough) of the economy during the introduction of innovations, the stage of slowing down and the exit of innovation from the economy \cite{Makarov_2014}.

There are several mathematical models of Kondratyev's cycles, the most common among them is the model of V.S. Dubovsky \cite{Dubovsky_1995} and the model of Akaev A.A. \cite{Akayev_2008}. We will focus on the model of V.S. Dubovsky, a generalization of which to the case of hereditarity was given by the authors in \cite{Makarov_2016}.

This work is a continuation of work \cite{Makarov_2016}. In the model equation, a function was introduced that characterizes the influx of new management decisions into innovation, and a more accurate numerical method for analyzing the proposed model was used, in contrast to the explicit non-local finite-difference scheme \cite{Parovik_2014}.

\section{Some reduction from the theory of fractional calculus}

Here we will consider the main definitions from the theory of fractional calculus, in more detail its aspects can be studied in the books \cite{Kilbas_2006,Oldham_1974,Miller_1993}.

\textbf{Definition 1.} Fractional Riemann-Liouville integral of order $\alpha$:
	\begin{equation}
		I_{0t}^\alpha x\left( \tau  \right) = \frac{1}{{\Gamma \left( \alpha 
				\right)}}\int\limits_0^t {\frac{{x\left( \tau  \right)d\tau }}{{{{\left( {t -
								\tau } \right)}^{1 - \alpha }}}}} ,\alpha  > 0,t > 0,
		\label{Makarov1}
\end{equation}	
here $\Gamma\left(.\right)$~ is the gamma function.

The operator (\ref{Makarov1}) has the following properties:
\[
I_{0t}^0x\left( \tau  \right) = x\left( t \right), I_{0t}^\alpha I_0^\beta
x\left( \tau  \right) = I_{0t}^{\alpha  + \beta }x\left( \tau  \right), I_{0t}^\alpha I_0^\beta x\left( \tau  \right) = I_0^\beta I_{0t}^\alpha x\left(\tau  \right).
\]	

\textbf{Definition 2.} The fractional Gerasimov-Caputo derivative of order $\alpha$ has the form:
	\begin{equation}
		\partial _{0t}^\alpha x\left( \tau  \right) = \left\{ \begin{array}{l}
			\dfrac{1}{{\Gamma \left( {m - \alpha } \right)}}\int\limits_0^t
			{\dfrac{{{x^{(m)}}\left( \tau  \right)d\tau }}{{{{\left( {t - \tau }
								\right)}^{\alpha  + 1 - m}}}},0 \le m - 1 < \alpha  < m,} \\
			\dfrac{{{d^m}x\left( t \right)}}{{d{t^m}}},m \in N.
		\end{array} \right.
		\label{Makarov2}
	\end{equation}		

The operator (\ref{Makarov2}) has the following properties\cite{Kilbas_2006}:
\[
\partial _{0t}^\alpha I_{0t}^\alpha x\left( \tau  \right) = x\left( t
\right),I_{0t}^\alpha \partial _{0t}^\alpha x\left( \tau  \right) -
\sum\limits_{k = 0}^{i - 1} {\frac{{{x^{\left( k \right)}}\left(0\right){t^k}}}{{k!}}}, t > 0.
\]

\section{Statement of the problem}

Consider the following Cauchy problem:

\begin{equation}
\left\{ \begin{array}{l}\partial _{0t}^{{\alpha _1}}x\left( \tau  \right) =  - \lambda nx\left( t\right)\left( {x\left( t \right) - 1} \right)\left( {y\left( t \right) - {y^ * }}
\right) + {\delta _1}\cos \left( {{\omega _1}t} \right),x\left( 0 \right) = a\\
\partial _{0t}^{{\alpha _2}}y\left( \tau  \right) = n\left( {1 - n}\right){y^2}\left( t \right)\left( {x\left( t \right) - {x^ * }} \right) +	{\delta _2}\cos \left( {{\omega _2}t} \right),y\left( 0 \right) = b\end{array} \right.,
\label{Makarov3} 
\end{equation}
where $ x \left (t \right) $ is the efficiency of innovation, the ratio of labor productivity at new jobs to average productivity at all jobs of all ages; $ y \left(t\right) $ -- efficiency of fixed assets (funds) organizations; $ n $ -- rate of accumulation, gross capital formation as a share of GDP; $ \lambda$ -- parameter that determines the size and duration of the cycles; $ a, b $ -- positive constants that determine the initial conditions; $ t \in \left [{0, T} \right] $ -- current time of the considered process; $ T> 0$ -- simulation time; $ \left ({{x^*}, {y^*}} \right)$ -- coordinate of the equilibrium point of system (\ref{Makarov1}); $ {\delta_1},{\delta_2}, {\omega _1}, {\omega _2}$ -- given positive constants; fractional operators in system (\ref{Makarov1}) are determined from (\ref{Makarov2}):
\begin{equation}
\partial _{0t}^{{\alpha _1}}x\left( \tau  \right) = \frac{1}{{\Gamma \left( {1 -{\alpha _1}} \right)}}\int\limits_0^t {\frac{{\dot x\left( \tau  \right)d\tau}}{{{{\left( {t - \tau } \right)}^{{\alpha _1}}}}}} ,0 < {\alpha _1} < 1, \partial _{0t}^{{\alpha _2}}y\left( \tau  \right) = \frac{1}{{\Gamma \left( {1 -{\alpha _2}} \right)}}\int\limits_0^t {\frac{{\dot y\left( \tau  \right)d\tau}}{{{{\left( {t - \tau } \right)}^{{\alpha _2}}}}}}, 0 < {\alpha _2} < 1.
\label{Makarov4}
\end{equation}

\textbf{Remark 1.} System (\ref{Makarov3}) describes Kondratyev's long waves taking into account heredity. In the case when in system (\ref{Makarov3}) $ {\alpha _1} = {\alpha _2} = 1 $ and $ {\delta _1} = {\delta _2} = 0 $ we obtain the Dubovskiy model, which was proposed and investigated S. V. Dubovsky in \cite{Dubovsky_1995}. Therefore, the dynamical system (\ref{Makarov3}) will be called the generalized Dubovsky model (GDM).

\section{Adams-Bashforth-Moulton method}

The ABM method is a type of numerical predictor-corrector method for solving differential equations. It has been studied and discussed in detail in the papers \cite{Garrappa_2018,Yang_2006,Diethelm_2002}. Let's generalize this method for solving the Cauchy problem (\ref{Makarov3}). 

We will assume that the required functions $ x \left (t \right), y \left (t \right) $ possess the required smoothness.  On a uniform grid, we introduce the functions $ x_{n + 1}^p, \, y_{n + 1}^p $, $n = 0, \ldots, \, N-1 $, which will be determined by the Adams-Bashforth formula (predictor):
\begin{equation}
\left\{ \begin{array}{l}
x_{n + 1}^p = {x_0} + \dfrac{{{\tau ^{{\alpha _1}}}}}{{\Gamma \left( {{\alpha _1}
			+ 1} \right)}}\sum\limits_{j = 0}^n {\theta _{j,n + 1}^1\left( { - \lambda
		n{x_j}\left( {{x_j} - 1} \right)\left( {{y_j} - {y^ * }} \right) + {\delta
			_2}\cos \left( {{\omega _2}j\tau } \right)} \right),} \\
y_{n + 1}^p = {y_0} + \dfrac{{{\tau ^{{\alpha _2}}}}}{{\Gamma \left( {{\alpha _2}
			+ 1} \right)}}\sum\limits_{j = 0}^n {\theta _{j,n + 1}^2\left( {n\left( {1 - n}
		\right)y_j^2\left( {{x_j} - {x^ * }} \right) + {\delta _2}\cos \left( {{\omega
				_2}j\tau } \right)} \right),} \\
\theta _{j,n + 1}^i = {\left( {n - j + 1} \right)^{{\alpha _i}}} - {\left( {n -
		j} \right)^{{\alpha _i}}}, i = 1,2.
\end{array} \right.
\label{Makarov5}
\end{equation}

Then, using the Adams-Moulton formula for the corrector, we get:
\begin{equation}
\left\{ \begin{array}{l}
{x_{n + 1}} = {x_0} + \dfrac{{{\tau ^{{\alpha _1}}}}}{{\Gamma \left( {{\alpha _1}
			+ 2} \right)}}\left( \begin{array}{l}
\left( { - \lambda nx_{n + 1}^p\left( {x_{n + 1}^p - 1} \right)\left( {y_{n +
			1}^p - {y^ * }} \right) + {\delta _1}\cos \left( {{\omega _1}\tau \left( {n + 1}
		\right)} \right)} \right) + \\
+ \sum\limits_{j = 0}^n {\rho _{j,n + 1}^1\left( { - \lambda n{x_j}\left(
		{{x_j} - 1} \right)\left( {{y_j} - {y^ * }} \right) + {\delta _1}\cos \left(
		{{\omega _1}j\tau } \right)} \right)}
\end{array} \right),\\
{y_{n + 1}} = {y_0} + \dfrac{{{\tau ^{{\alpha _2}}}}}{{\Gamma \left( {{\alpha _2}
			+ 2} \right)}}\left( \begin{array}{l}
\left( {n\left( {1 - n} \right){{\left( {y_{n + 1}^p} \right)}^2}\left( {x_{n +
			1}^p - {x^ * }} \right) + {\delta _2}\cos \left( {{\omega _2}\tau \left( {n + 1}
		\right)} \right)} \right) + \\
+ \sum\limits_{j = 0}^n {\rho _{j,n + 1}^2\left( {n\left( {1 - n}
		\right)y_j^2\left( {{x_j} - {x^ * }} \right) + {\delta _2}\cos \left( {{\omega
				_2}j\tau } \right)} \right)}
\end{array} \right),
\end{array} \right.
\label{Makarov7}
\end{equation}
where the weight factors in (\ref{Makarov7}) are determined by the formula:
\[
\rho _{j,n + 1}^i = \left\{ \begin{array}{l}
{n^{{\alpha _i} + 1}} - \left( {n - {\alpha _i}} \right){\left( {n + 1}\right)^{{\alpha _i}}},\,j = 0,\\ {\left( {n - j + 2} \right)^{{\alpha_i} + 1}} + {\left( {n - j}
	\right)^{{\alpha_i} + 1}} - 2{\left( {n - j + 1} \right)^{{\alpha _i} + 1}},\,1\le j \le n,\\ 1,\,j = n + 1,\\ i = 1,2. \end{array} \right.
\]

\textbf{Theorem \cite{Yang_2006}.} If $\partial _{0t}^{{\alpha_i}}{x_i}\left( \tau  \right) \in {C^2}\left[ {0,T} \right],\left( {{x_1} = x\left( t \right),\,{x_2} = y\left(t\right),\,i = 1,2} \right)$, then
	\begin{equation}
	\mathop {\max }\limits_{1 \le j \le n} \left|{{x_i}\left( {{t_j}} \right) - {x_{i,j}}} \right| = O\left( {{h^{1 + \mathop
				{\min }\limits_i {\alpha_i}}}} \right).
	\label{Makarov8}
	\end{equation}
	
The proof of the theorem is based on the method of mathematical induction, and it is given in \cite{Yang_2006}.

\textbf{Remark 2.} Note that in the case $ {\alpha _i} = 1 $, taking into account (\ref{Makarov7}), we obtain the classical ABM method of the second order of accuracy.

\textbf{Remark 3.} In \cite{Makarov_2016}, an explicit finite difference scheme was used for research, which has conditional convergence and stability. The ABM method scheme is devoid of these disadvantages.	
	
\section{Research results}

Let us examine how the computational accuracy of the methods behaves. To do this, we will use the double recalculation method (Runge's rule) to estimate the error using the formula:
\begin{equation}
\xi_x^i = \frac{{\mathop {\max }\limits_i \left|{{x_i} - {x_{2i}}} \right|}}{{{2^\mu } - 1}},\xi _y^i = \frac{{\mathop {\max}\limits_i \left| {{y_i} - {y_{2i}}} \right|}}{{{2^\mu } - 1}}, 
\end{equation}
where $ \mu = 1 + \min \left ({{\alpha _1}, {\alpha _2}} \right) $ is the order of accuracy of the ABM method, $ \xi _x^i, \xi_y^i $ are errors at step $ {\tau _i} $, $ {x_ {2i}}, {y_ {2i}} $ -- numerical solutions at step $ {{{\tau _i}} \mathord {\left / {\vphantom {{{\tau _i}} 2}} \right. \kern \nulldelimiterspace} 2} $.

 The computational accuracy of $ p $ is determined from the formula:
\begin{equation}
{p_x} = {\log _2}\left( {{{\xi _x^i}\mathord{\left/{\vphantom {{\xi _x^i} {\xi _x^{i + 1}}}} \right.\kern-\nulldelimiterspace} {\xi _x^{i + 1}}}} \right),    {p_y} = {\log_2}\left( {{{\xi _y^i} \mathord{\left/{\vphantom {{\xi _y^i} {\xi _y^{i + 1}}}} \right.\kern-\nulldelimiterspace} {\xi _y^{i + 1}}}} \right),
\label{Makarov9}		                             
\end{equation}
$ {\tau_i}, {\tau_ {i + 1}} = {{{\tau _i}} \mathord {\left/{\vphantom {{{\tau _i}} 2}} \right. \kern-\nulldelimiterspace} 2} $ are grid steps, and $ \xi_x^{i + 1}, \xi_y^{i + 1} $ at step $ {\tau_{i + 1}}$.

\textbf{Example 1.} (Classical Dubovsky model).

Let's consider a test case in the case $ {\alpha _1} = {\alpha _2} = 1 $. The values of the parameters for the Cauchy problem (\ref{Makarov3}) are chosen as follows: $ {\delta _1} = {\delta _2} = 1, {\omega _1} = {\omega _2} = 0.5, n = 0.2, \lambda = 1.5, x \left (0 \right) = 5, y \left (0 \right) = 4, {x^*} = 1.35, {y ^ *} = 0.5, t \in \left [{0,1} \right ] $. The results of the numerical analysis are shown in Table 1.

\begin{table}[h!]
	\caption{The classical Dubovsky model $ \left ({{\alpha _1} = {\alpha _2} = 1} \right) $}
	\begin{ruledtabular}
	\begin{tabular}{|c|c|c|c|c|c|}
		$N$ &  $\tau $ & ${\xi_x^i}$ & ${\xi_y^i}$ &  ${p_x}$ & ${p_y}$ \\
		\hline
		
		     10
		 &  
			1/10
		 &  
			0.0397415350
		 &  
			0.0799740787
		 &  
			-
		 &  
			-
		 \\
		 
			20
		 & 
			1/20
		 & 
			0.0060318953
		 & 
			0.0187709253
		 &  
			2.72
		 & 
			2.09
      	 \\
		
			40
		 &  
			1/40
		 &  
			0.0010575820
		 &  
			0.0045870660
		 &  
			2.51
		 & 
			2.03
		 \\
		 
			80
		 &  
			1/80
		 & 
			0.0002336673
		 &  
			0.0011345727
		 & 
			2.18
		 & 
			2.02
		 \\
		
			160
		 &  
			1/160
		 & 
			0.0000552493
		 &  
			0.0002821520
		 &  
			2.08
		 &  
			2.01
		 \\
		
			320
		 & 
			1/320
		 &  
			0.000013453
		 & 
			0.0000703530
		 & 
			2.04
		 &  
			2.00
		 \\
	\end{tabular}
	\end{ruledtabular}
\end{table}

Table 1 shows that with an increase in the calculated grid nodes, the computational accuracy $ {p_x}, {p_y} $ tend to $ \mu = 2 $ is the order of accuracy of the ABM method. Consider another example that implements ODM.

\textbf{Example 2.} (Generalized Dubovsky model)

Consider the case $ {\alpha _1} = 0.9, {\alpha _2} = 0.8 $, the rest of the parameters are taken from the previous example. The simulation results are shown in Table 2.

\begin{table}[h!]
\caption{Generalized Dubovsky model $ \left ({{\alpha _1} = 0.9, {\alpha _2} = 0.8} \right)$}
\begin{ruledtabular}
	\begin{tabular}{|p{17pt}|p{30pt}|p{77pt}|p{77pt}|p{77pt}|p{77pt}|}
		\parbox{17pt}{\centering $N$} & \parbox{30pt}{\centering $\tau $} & \parbox{77pt}{\centering ${\xi_x^i}$} & \parbox{77pt}{\centering ${\xi_y^i}$} & \parbox{77pt}{\centering ${p_x}$} & \parbox{77pt}{\centering ${p_y}$} \\
		\hline
		\parbox{17pt}{\centering 
			10
		} & \parbox{30pt}{\centering 
			1/10
		} & \parbox{77pt}{\centering 
			0.1171054420
		} & \parbox{77pt}{\centering 
			0.2158809510
		} & \parbox{77pt}{\centering 
			-
		} & \parbox{77pt}{\centering 
			-
		} \\
		\parbox{17pt}{\centering 
			20
		} & \parbox{30pt}{\centering 
			1/20
		} & \parbox{77pt}{\centering 
			0.0235086005
		} & \parbox{77pt}{\centering 
			0.0586392164
		} & \parbox{77pt}{\centering 
			2.3165475629
		} & \parbox{77pt}{\centering 
			1.8802982187
		} \\
		\parbox{17pt}{\centering 
			40
		} & \parbox{30pt}{\centering 
			1/40
		} & \parbox{77pt}{\centering 
			0.0041498810
		} & \parbox{77pt}{\centering 
			0.0173332737
		} & \parbox{77pt}{\centering 
			2.5020467762
		} & \parbox{77pt}{\centering 
			1.7583216658
		} \\
		\parbox{17pt}{\centering 
			80
		} & \parbox{30pt}{\centering 
			1/80
		} & \parbox{77pt}{\centering 
			0.0009959756
		} & \parbox{77pt}{\centering 
			0.0051032771
		} & \parbox{77pt}{\centering 
			2.0588875967
		} & \parbox{77pt}{\centering 
			1.7640482636
		} \\
		\parbox{17pt}{\centering 
			160
		} & \parbox{30pt}{\centering 
			1/160
		} & \parbox{77pt}{\centering 
			0.0002742665
		} & \parbox{77pt}{\centering 
			0.0014889669
		} & \parbox{77pt}{\centering 
			1.8605318928
		} & \parbox{77pt}{\centering 
			1.7771123047
		} \\
		\parbox{17pt}{\centering 
			320
		} & \parbox{30pt}{\centering 
			1/320
		} & \parbox{77pt}{\centering 
			0.0000779949
		} & \parbox{77pt}{\centering 
			0.0004308895
		} & \parbox{77pt}{\centering 
			1.8141277492
		} & \parbox{77pt}{\centering 
			1.7889216854
		} \\
	\end{tabular}
\end{ruledtabular}	
\end{table}

From Table 2 we see that the computational accuracy tends to the order of accuracy of the ABM method, which has a value of $ \mu = 1.8 $.

The generalized mathematical model of Dubovsky (\ref{Makarov3}) has both regular and chaotic modes. Let us show this using the numerical algorithm of the ABM method and construct phase trajectories for different values of the problem parameters.

Figure 1 shows the regular modes of the generalized mathematical model of Dubovsky. Figure 1a shows the phase trajectories constructed depending on different values of $ \lambda $ and initial conditions $ \left ({{\delta_1} = {\delta _2} = 0} \right ) $, which were taken from \cite{Dubovsky_1995}, that is, the classical Dubovsky model $ {\alpha _1} = {\alpha _2} = 1 $ is given.

\begin{figure}[h!]
\begin{minipage}[h]{0.24\linewidth}
{\includegraphics[scale=0.22]{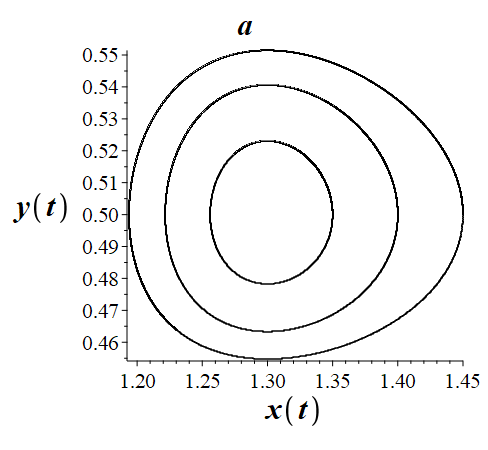}}
\end{minipage}
\begin{minipage}[h]{0.24\linewidth}
{\includegraphics[scale=0.22]{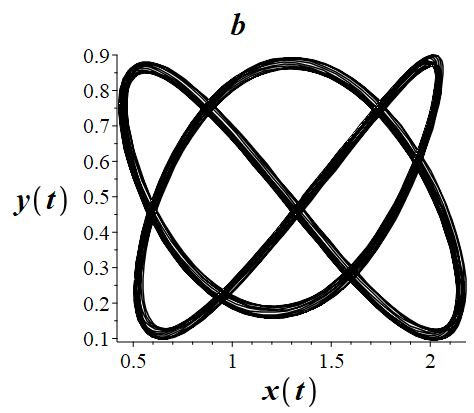}}\\
\end{minipage}
\vfill
\begin{minipage}[h]{0.24\linewidth}
{\includegraphics[scale=0.22]{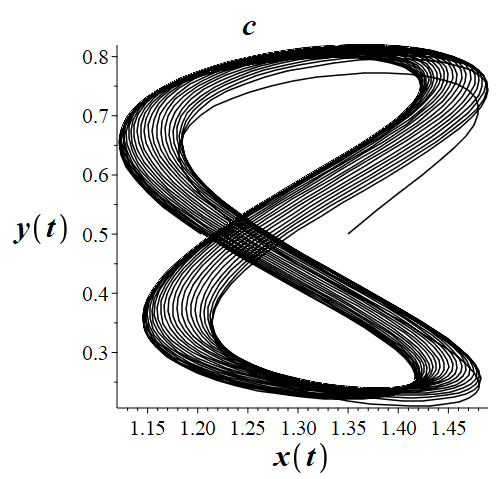}}
\end{minipage}
\begin{minipage}[h]{0.24\linewidth}
{\includegraphics[scale=0.22]{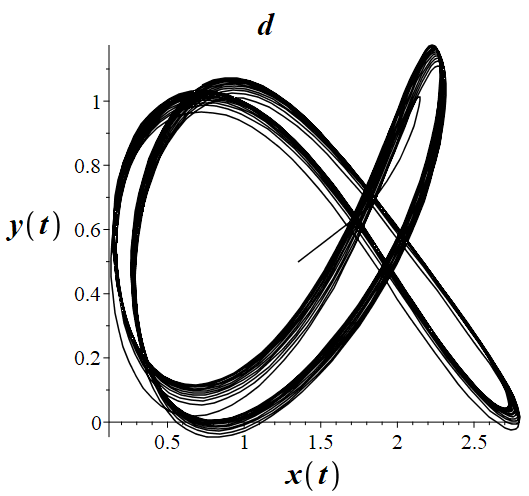}}
\end{minipage}
\caption{Regular modes}
\end{figure}

It can be seen that the phase trajectories in Figure 1a are in good agreement with the results of \cite{Dubovsky_1995}, which indicates the correctness of calculations using the ABM method. Figure 1b, c, d shows the phase trajectories obtained for various values of $ {\alpha_1}, {\alpha _2}, {\delta _1}, {\delta _2}, {\omega _1}, {\omega _2} , \lambda $. We see that the phase trajectories are limit cycles, but of a more complex shape than in the classical case.

\begin{figure}[h!]
\begin{minipage}[h]{0.3\linewidth}
	{\includegraphics[scale=0.24]{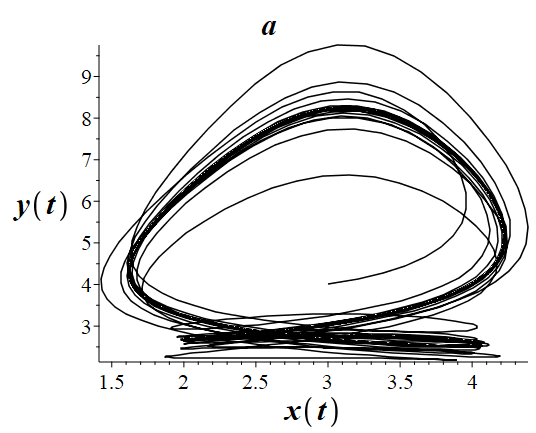}}
\end{minipage}
\begin{minipage}[h]{0.3\linewidth}
	{\includegraphics[scale=0.24]{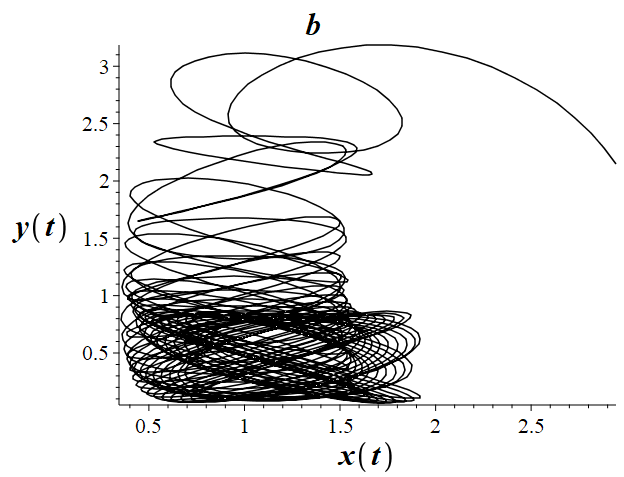}}
\end{minipage}
\caption{Chaotic modes}
\end{figure}

Figure 2 shows the phase trajectories that characterize chaotic regimes. Figure 2a is the origin of a chaotic attractor, Figure 2b is a chaotic attractor. Such modes need to be investigated in more detail, for example, by analogy with works \cite{Parovik_2018,Cao_2020,Diouf_2020}. Such a study will allow taking into account the necessary values of the model parameters for the existence of limit cycles.

\section{Conclusion}

In this paper, we have proposed a generalized Dubovsky model, taking into account the effects of heredity, as well as functions responsible for investment and management technologies. We examined the ABM method, investigated the accuracy of the method, and built phase trajectories. They showed that solutions can describe both regular regimes and chaotic regimes.

Of interest is the further study of the model in the following areas: the study of chaotic regimes, for example, using the spectra of the maximum Lyapunov exponents, the study of limit cycles, the determination of their lengths, for example, using the Poincaré sections, the economic interpretation of research results.

\begin{acknowledgments}
The work was performed within the framework of the research project of Vitus Bering KamSU "Mathematical model of Kondratiev's long waves taking into account heredity" No. AAAA-A20-120021190003-1.
\end{acknowledgments}

\nocite{*}
\bibliography{Makarov}

\end{document}